\documentclass[]{interact}

\usepackage{epstopdf}
\usepackage[caption=false]{subfig}

\usepackage[longnamesfirst,sort]{natbib}
\bibpunct[, ]{(}{)}{;}{a}{,}{,}
\renewcommand\bibfont{\fontsize{10}{12}\selectfont}

\theoremstyle{plain}
\newtheorem{theorem}{Theorem}[section]

\theoremstyle{definition}

\theoremstyle{remark}
\newtheorem{remark}{Remark}

\begin{document}


\title{Fold-and-cut theorem: a vertical maths project around origami}

\author{
\name{Maria Luisa Spreafico\textsuperscript{a}\thanks{CONTACT M.L. Spreafico. Email: maria.spreafico@polito.it} and Eulalia Tramuns\textsuperscript{a}}
\affil{\textsuperscript{a} Dipartimento di Scienze Matematiche, Politecnico di Torino, C.so Duca degli Abruzzi 24, 10129 Torino, Italy}
}

\maketitle

\begin{abstract}
Fold-and-cut theorem claims the possibility to cut out from a sheet a set of straight-line drawing using only one cut of scissors, without producing any other cut in the sheet and separating all the figures at the same time, just by folding flat the paper before the last cut. This fascinating and magic result, proved by E.Demaine and alii, allows a wide exploration of mathematical issues.
We present in this paper several different approaches about the use of this result, as problem based learning, direct investigation, manipulative exploration and use of ICT, most of them based on learning by doing strategy. Most of the lessons we present have been tested on a group fo 170 students of primary and middle school and around 30 adults, who showed to be mostly satisfied with the experience.

\end{abstract}


\section{Introduction}

The benefits of origami in learning mathematics are well-known. Touching models improves spatial visualization and is clearly helpful in geometric contexts, improving the implicit understanding with manipulative strategies, but also in algebra and calculus and other areas, as show several studies, see \cite{Arc15, Lam07, Rob03, Yuz02, Mey99}.

Moreover some works show practical examples of application of origami and describe varied lessons that can be easily reproduced in classes, for different levels of the educational system: some related with 3D geometry, areas and volumes, see \cite{War11}, with logic, see \cite{Ser18}, and also cross-curricular lessons involving for example art, mathematics and origami, see \cite{Spr19, Spr18}. 

Moreover, the use of origami in mathematics meets the multiple intelligences of the students, identified by Howard Gardner, see \cite{Gar06}.

But paper-folding is also useful to bring magic, suspense and enthusiasm. In this paper we describe how using this facet allows to propose mathematics from another point of view, detailing lessons focused on fold-and-cut theorem. In this context, we allow students to use scissors while it is not contemplate  in the traditional origami. This project gives an example of an interesting topic that can be proposed vertically in a school context involving pupils in relation to their competences.
Specifically, the participants of this project around maths and origami activities were students of 6 classes of primary school, that is roughly 150 students aged between 8 and 10, 27 students of middle school aged 12 of an Italian school.

Each class did 4 math lessons: 3 of them following the curricular program and one lesson out of it, according to the students' skills. We focus on this lesson because is an interesting example of how the same mathematical result can be addressed from different perspectives and levels, adapting the depth of the explanations and inserting the use of ICT when possible. 

The approaches we proposed where different according to the age of students. In the next section we describe in detail the activities. Table \ref{summary_table} collects the didactical approaches of the activities.

\begin{table}
\tbl{Activities and type of approach.}
{\begin{tabular}{lcccccc} \toprule
& \multicolumn{2}{l}{Activities} \\ \cmidrule{2-7}
Approach & One & Two & Three & Four & Five & Six \\ \midrule
Free exploration & X & &  & & & X \\
Manipulative exploration & X &  & X & X & X &  \\ 
Problem solving &  &  & X & X & X&  \\
Use of ICT &  & X &  &  &  & X  \\ \bottomrule
\end{tabular}}
\label{summary_table}
\end{table}

The mathematical objectives pursued by this project were: to use geometry knowledge to check a result and find a solution of it; to work in increasing steps of difficulty; to propose first formal idea of proof, while solving questions as
 \emph{is it possible to fold and cut a given polygon from a paper?}.

Apart from the mathematical objectives, we wanted to approach a modern math result to show that mathematics are a living science, and current results at university level can be explored and divulgated to young people, adapting the language to them.
All the photos of the paper are by the authors.

\section{Material and methods}

One-cut results are known since times. A proof of this is the famous magician Houdini who used it as a trick around 1900s. But the exploration of what encompasses these isolated examples is relatively recent and affirms the possibility to generalise the result in any straight-line drawing. 

Technically, the theorem claims (see \cite{Rou11}, Theorem 5.1):
\begin{theorem}[Fold-and-Cut]
Any straight-line drawing (one composed of straight segments) on a sheet of paper may be folded flat so that one straight scissors cut completely through the folding cuts all the segments of the drawing and nothing else. 
\end{theorem}

In fact, this theorem allows different kinds of problems, because there are two different sets of straight-lines to consider: the initial straight-line drawing itself and the set of folds that allow that the paper remains flat, that become creases on the paper. From a didactical point of view,  one can begin considering the straight-line drawing or think about an inverse problem, beginning from the creases and get/explore the initial set of straight-lines.

In the following we propose several examples of these applications, developped in 6 activities. 

\begin{itemize}
\item Activity 1: Free cut! 

\begin{figure}
\centering
\subfloat{%
\resizebox*{6.9cm}{!}{\includegraphics{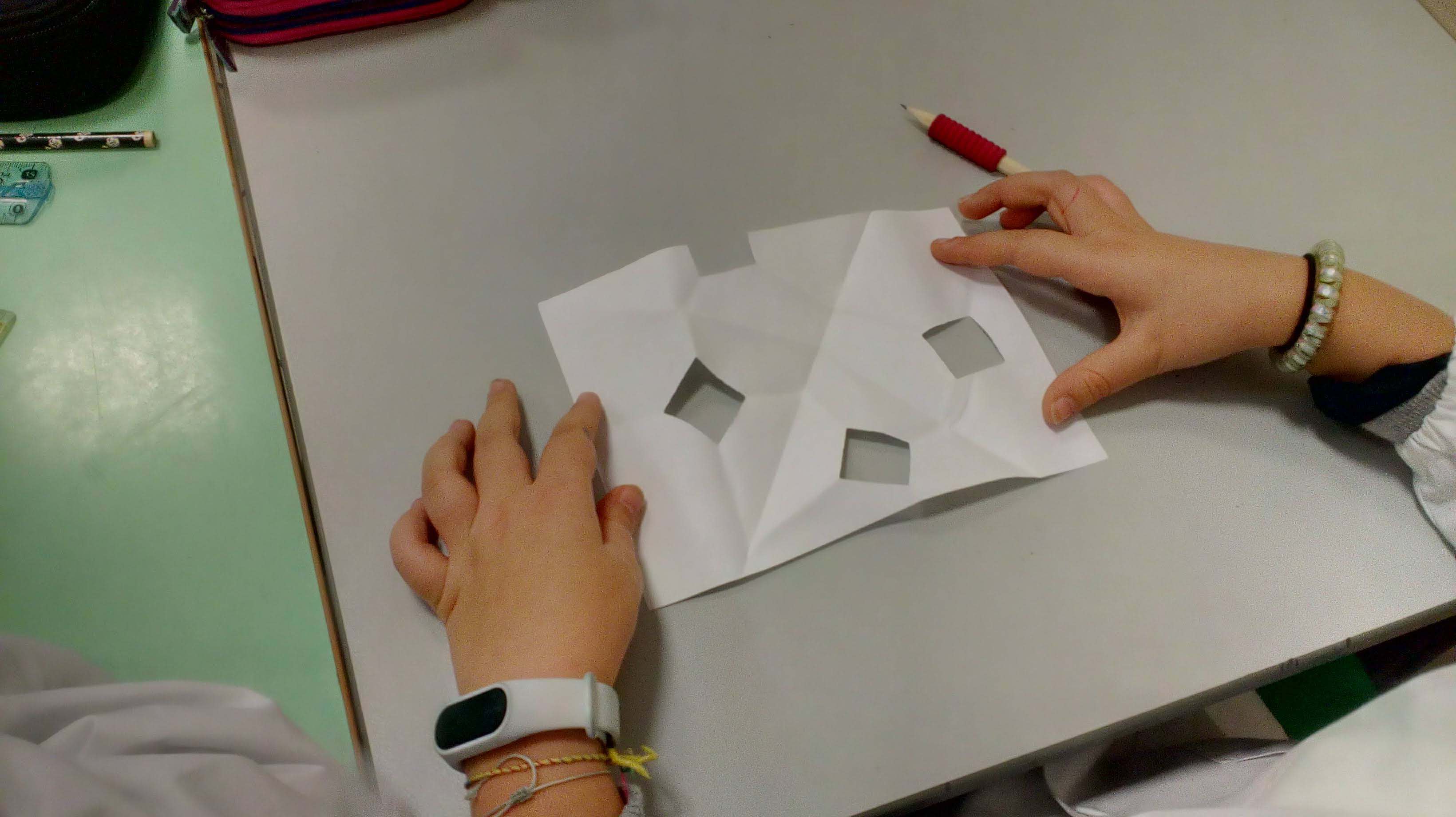}}}\hspace{5pt}
\subfloat{%
\resizebox*{6.3cm}{!}{\includegraphics{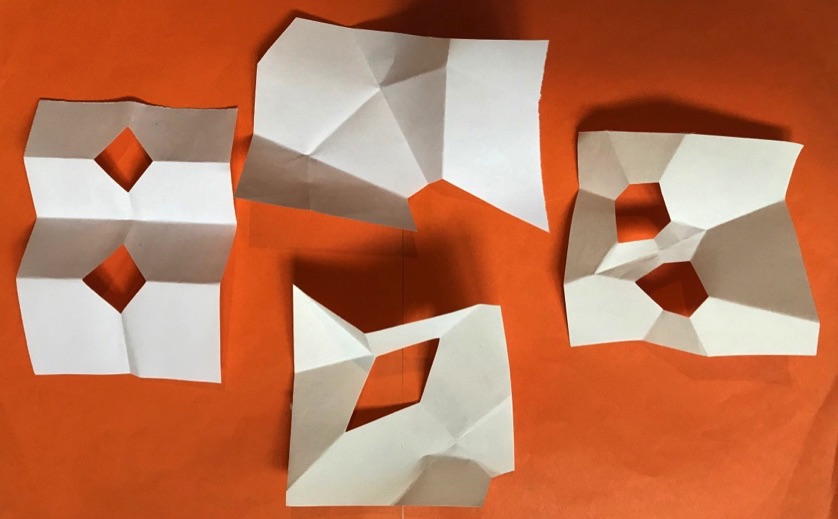}}}
\caption{Examples of free cut} \label{fig:FreeCut}
\end{figure}

Distribute a sheet of paper to each student. Let each student fold it flat freely. Each student should then draw a segment that starts and arrives at the edges of the paper and cut along it. Everyone will reopen the sheet and observe the different cuts produced. In some cases, there will be holes with symmetries, in others cuts on the outermost parts of the paper, see Figure \ref{fig:FreeCut} for some examples.

Free experimentation is important to let creativity and questions flow.

\item Activity 2: Video watching.  

Show some part of the interesting lesson by Prof. Demaine (see \cite{Dem12}) and claim the theorem. In this lesson he shows some surprising cases of Fold-and-cut theorem and children can understand that using mathematics it is possible to design a paper-folding pattern to get the required figures (such as a turtle or a butterfly). 

Considering the video and examples on Prof. Demaine site, one has to be aware of some origami vocabulary. I.e. linear folds are separed in two types, mountain and valley folds, see Figure \ref{fig:MountainValley}(a) and (b). The set of marks is named crease pattern. The notation for mountain folds is dash-dotted and for valley dashed, see Figure \ref{fig:MountainValley}(c).

\begin{figure}
\centering
\subfloat[Mountain fold.]{%
\resizebox*{3.1cm}{!}{\includegraphics{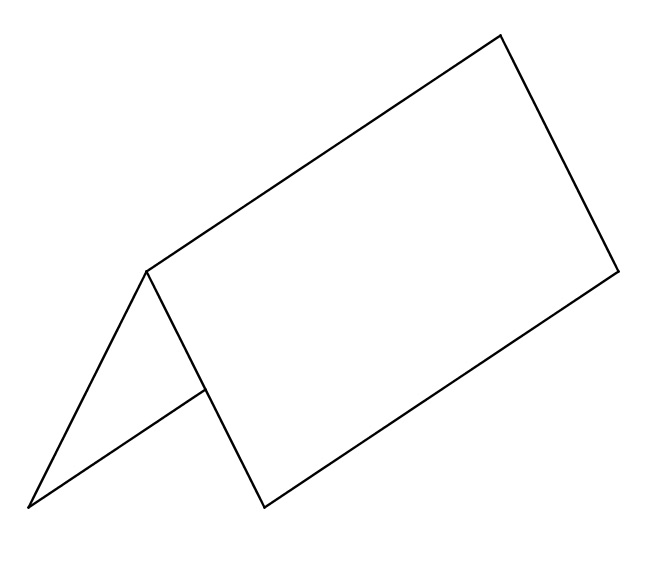}}}\hspace{5pt}
\subfloat[Valley fold.]{%
\resizebox*{3.1cm}{!}{\includegraphics{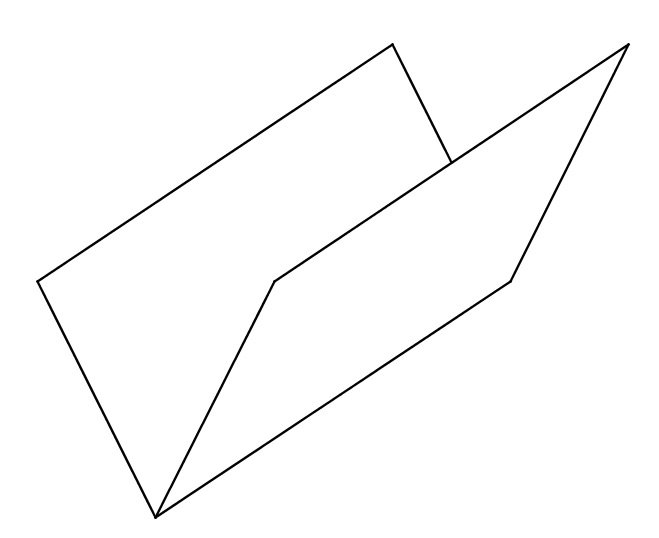}}}\hspace{5pt}
\subfloat[Crease pattern (CP)]{%
\resizebox*{3.1cm}{!}{\includegraphics{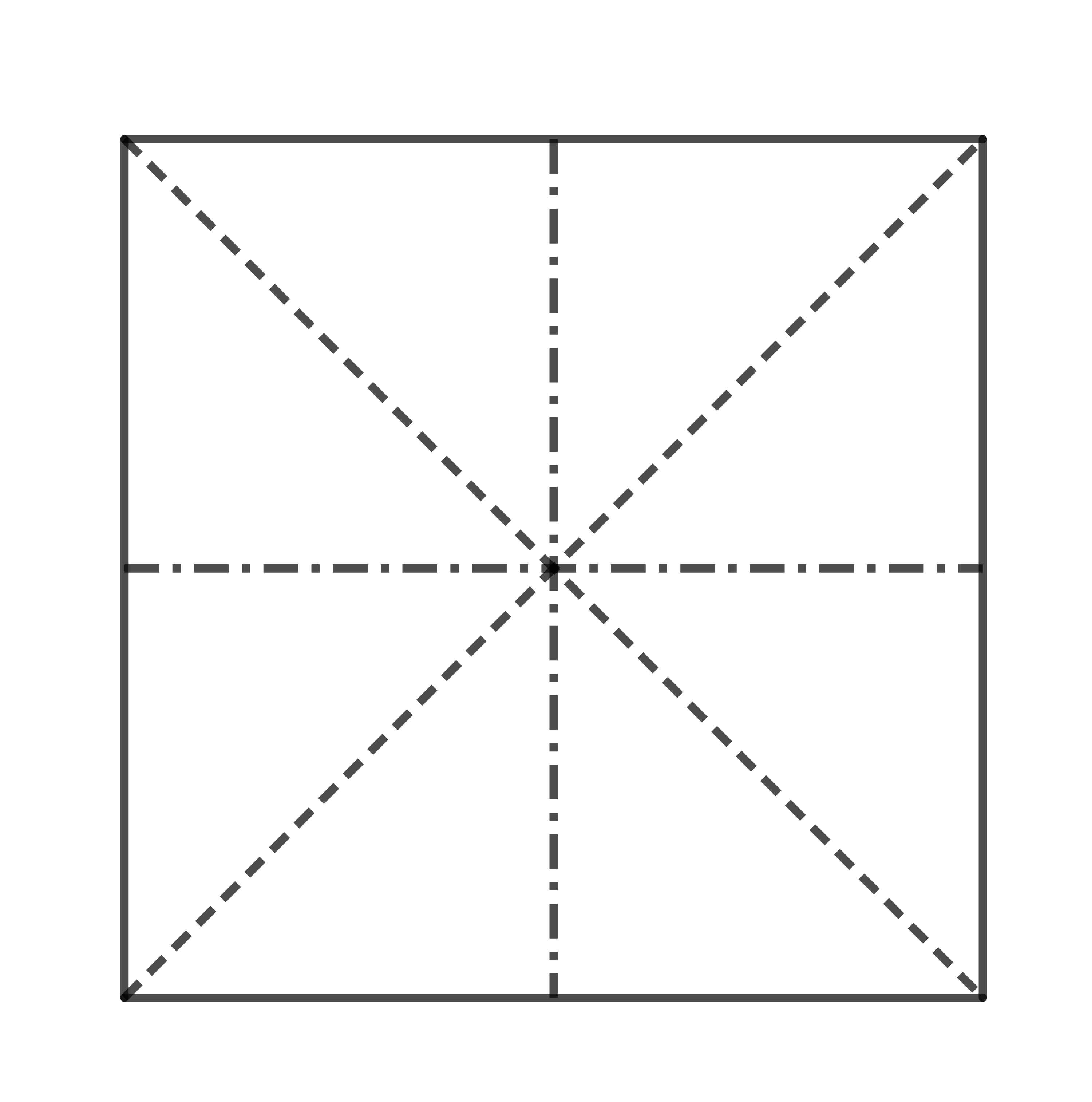}}}
\caption{Origami terminology} \label{fig:MountainValley}
\end{figure}

\item Activity 3: Design a way to cut out a square from a sheet. 

In this simple case children have to experiment to find the folding steps that allow to separate a square drawn on paper from the sheet. They have to understand that the use of symmetry axes help them to do it. From our experience, most of them are able to find at least two different folding sequences, see Figure \ref{fig:Solutions}:
\begin{figure}
\centering
\subfloat{%
\resizebox*{6.3cm}{!}{\includegraphics{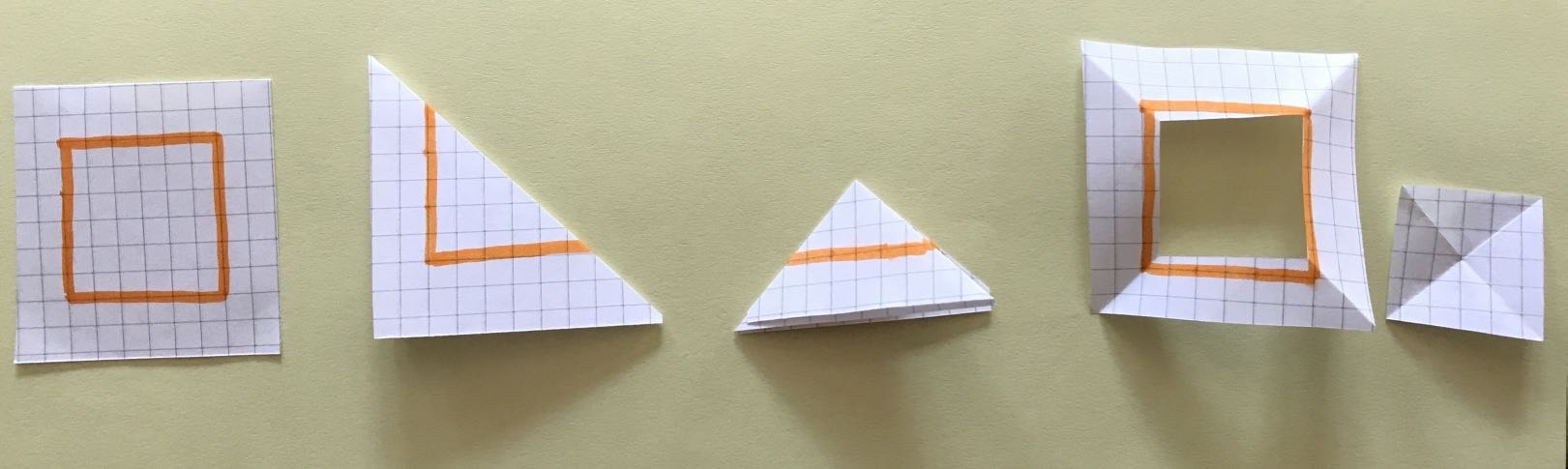}}}\hspace{5pt}
\subfloat{%
\resizebox*{6.8cm}{!}{\includegraphics{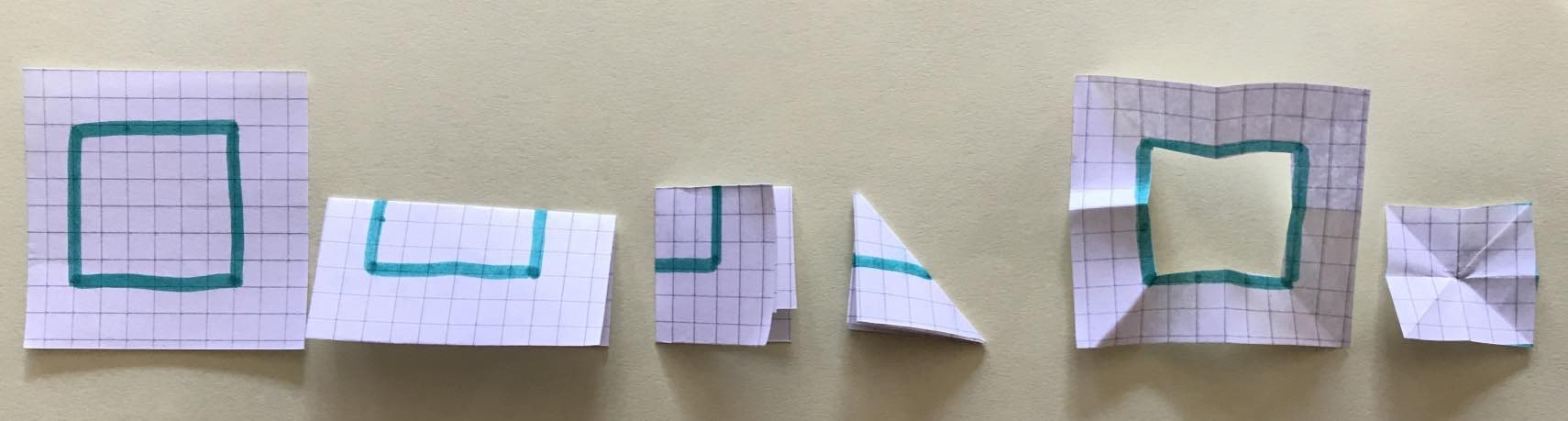}}}
\caption{Two different ways to cut a square} \label{fig:Solutions}
\end{figure}

\begin{itemize}
\item[(i)] they fold the sheet along the two diagonals of the square, one after the other, without unfolding the first, and cut the square side (2 folds and 1 cut),  
\item[(ii)] they fold a symmetry axe parallel to a side of the square, obtaining a rectangle, then they fold its shorter axe of symmetry, obtaining a square, and finally they fold a diagonal of this last small square cutting its side (3 folds 1 cut).
\end{itemize}

Figure \ref{fig:Trial1} show examples of trials of primary school students with no success. These cases are very interesting and helped to identify that some concepts were not assimilated, such that symmetry axes. Figure \ref{fig:Success}(a) shows the trail of another student in the same class, who succeeded.

\begin{figure}
\centering
\subfloat{%
\resizebox*{7cm}{!}{\includegraphics{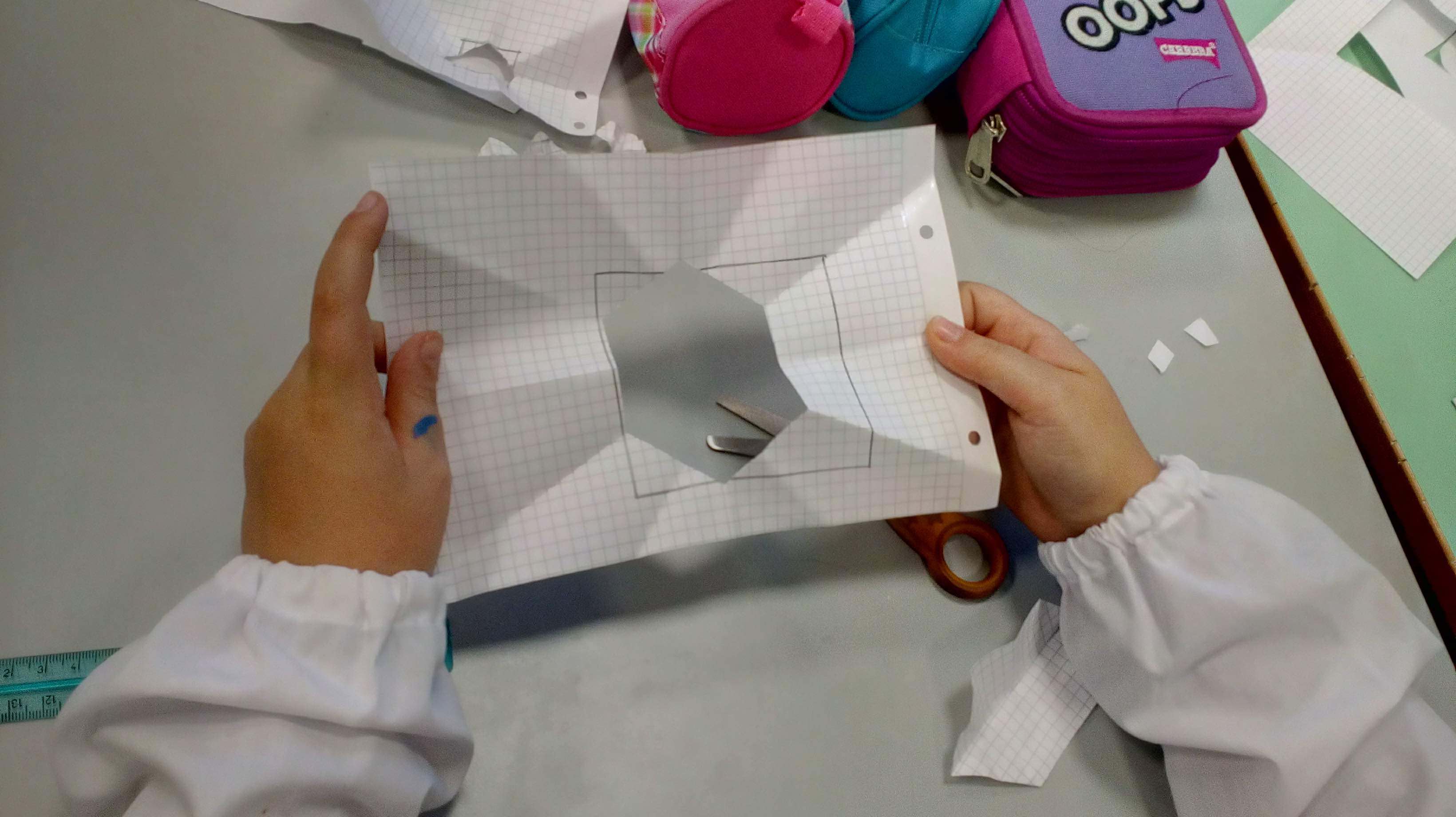}}}\hspace{5pt}
\subfloat{%
\resizebox*{7cm}{!}{\includegraphics{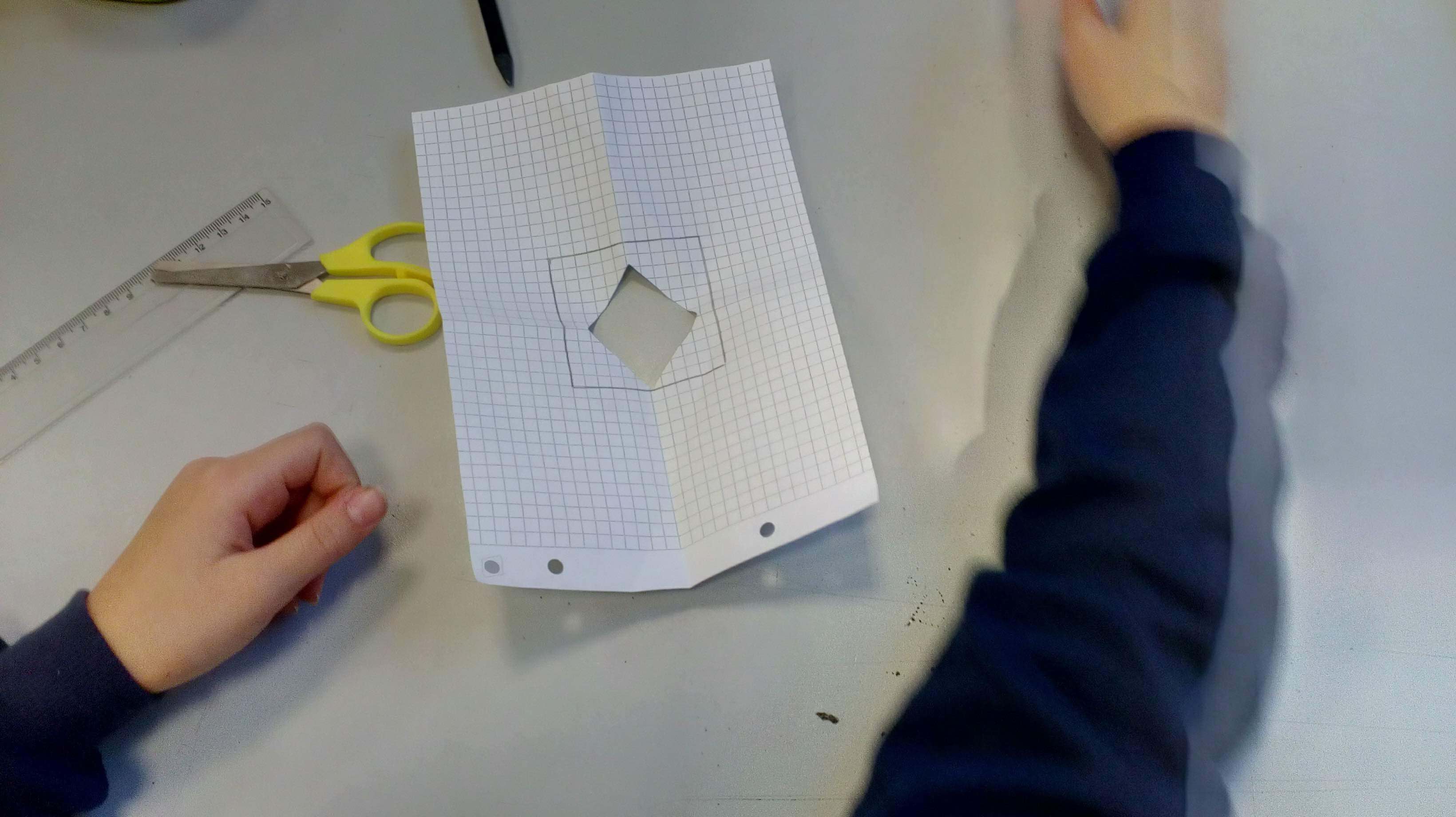}}}
\caption{Examples of trials to cut the square with no success} \label{fig:Trial1}
\end{figure}


\begin{figure}
\centering
\subfloat[]{%
\resizebox*{7cm}{!}{\includegraphics{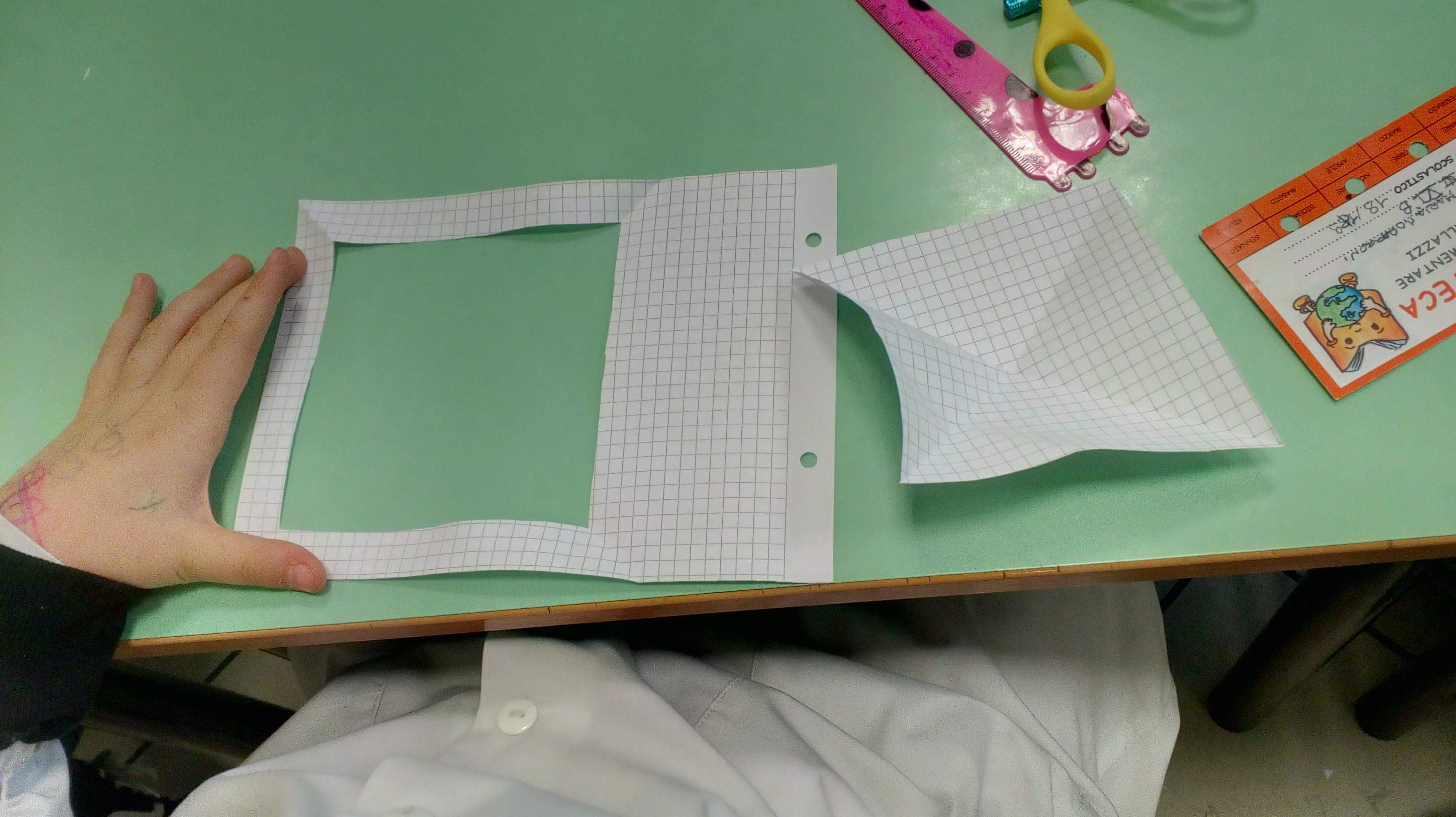}}}\hspace{5pt}
\subfloat[]{%
\resizebox*{5.25cm}{!}{\includegraphics{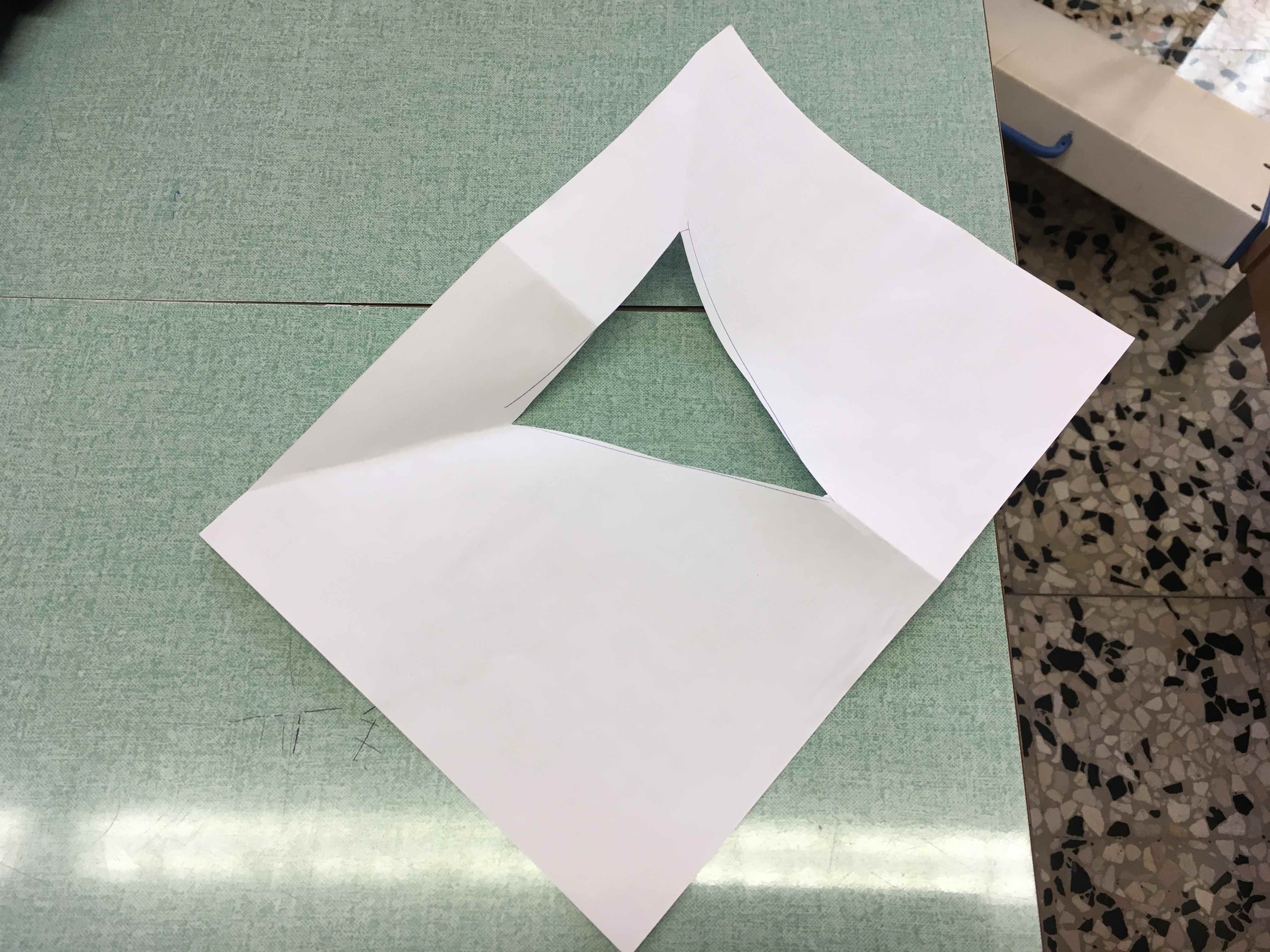}}}
\caption{Examples of trials to cut a square and a triangle} \label{fig:Success}
\end{figure}


This easy exercise, with its different possibilities, allows the teacher to speak about the variety and the convenience of solutions in maths.

\item Activity 4: Funny homework! 

We propose to kids, as experience at home, to fold-and-cut out in one fell swoop an ordinary triangle and a rectangle, one each time, explaining their choice during the folding process. 
A solution got by a 12 years old student is in Figure \ref{fig:Success}(b), when the triangle is isosceles.


\begin{figure}
\centering
\subfloat[]{%
\resizebox*{7.2cm}{!}{\includegraphics{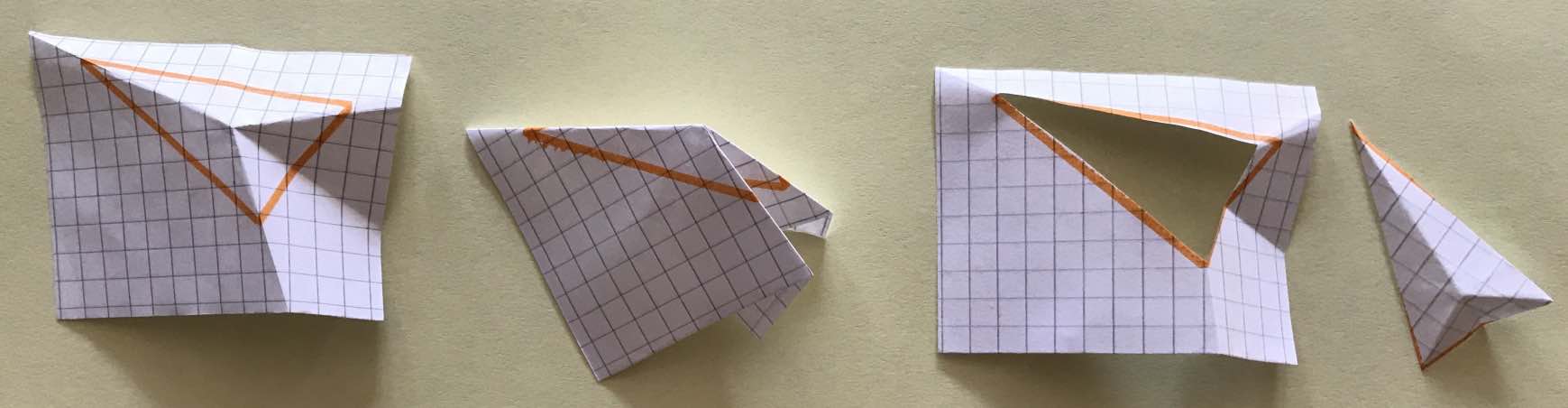}}}\hspace{5pt}
\subfloat[]{%
\resizebox*{5.5cm}{!}{\includegraphics{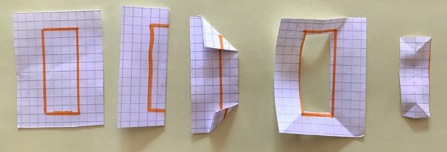}}}
\caption{Folding and cutting a triangle and a rectangle} \label{fig:TriangleRectangle}
\end{figure}

One possible solution for cutting out the triangle is to fold the three bisectors, see Figure \ref{fig:TriangleRectangle}(a). Regarding the rectangle one option is to fold upstream along the longest median and then fold the two halves of the short sides on the longest visible side. Then cut on the segment given by the long side, see Figure \ref{fig:TriangleRectangle}(b). Note how the bisectors of the vertex corners (bent upstream in the figure) plus a segment perpendicular to one side and passing through the incenter (downstream in the figure) are used.

\item Activity 5. Fold the diagram and cut: what does it appear?


\begin{figure}
\centering
\subfloat[]{%
\resizebox*{5.7cm}{!}{\includegraphics{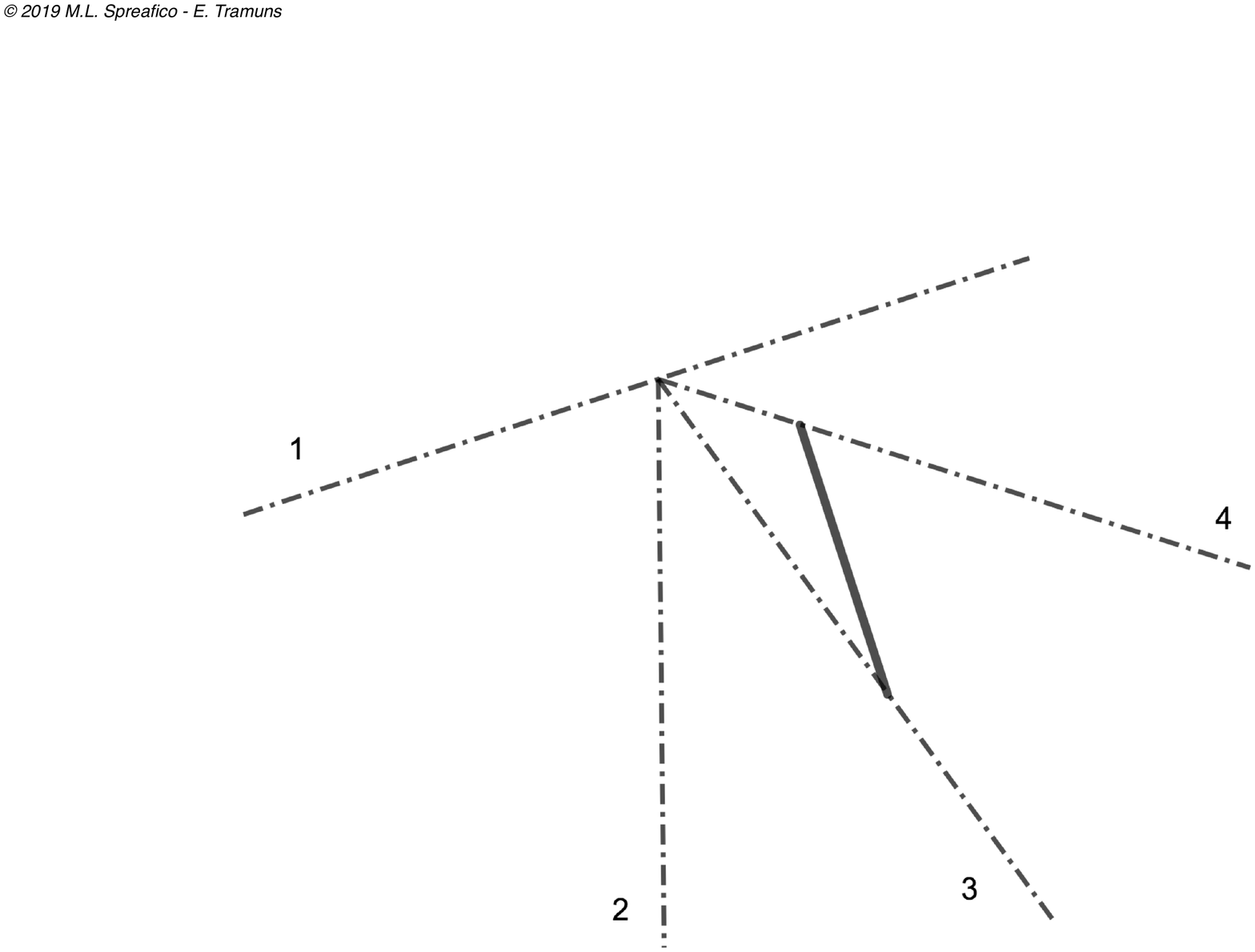}}}\hspace{5pt}
\subfloat[]{%
\resizebox*{5.7cm}{!}{\includegraphics{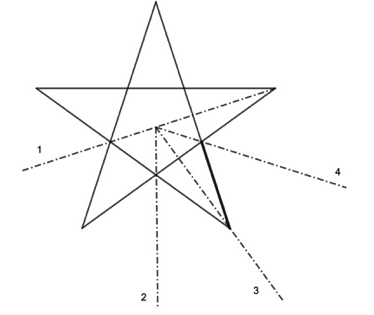}}}
\caption{The star pattern} \label{fig:starpattern}
\end{figure}

\begin{figure}[h!tbp]
\centering
	\includegraphics[width=0.7\textwidth]{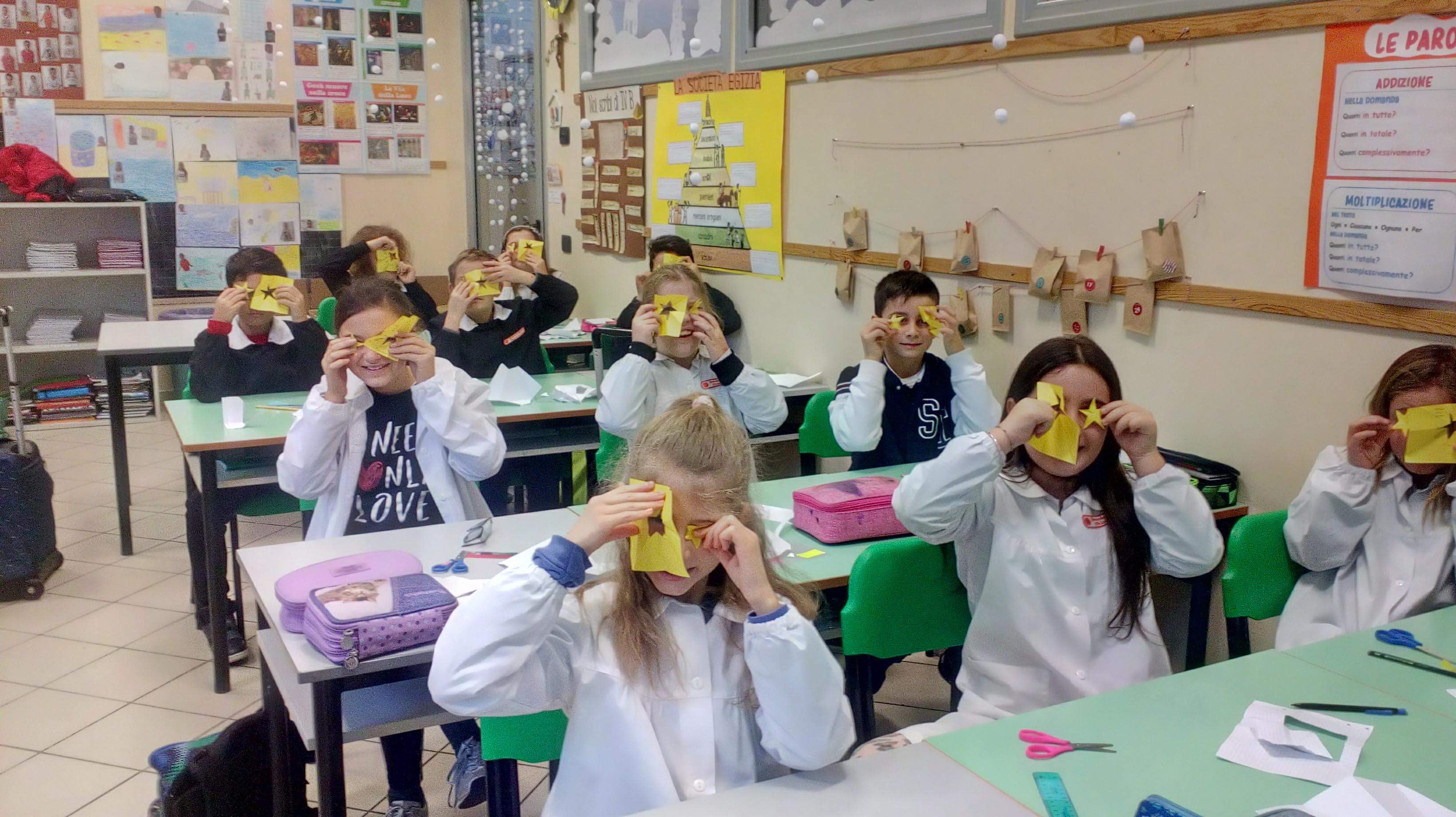}
\caption{Cutting stars; primary school level} \label{fig:Starprimary}
\end{figure}

In this last part we want to surprise children. We distribute papers with creases, see Figure \ref{fig:starpattern}(a).
Then the teacher gives instructions: everyone has to follow the order given by the numbers, taking into account that all the folds are mountain folds (dashed-dotted folds), that means that we have to fold to the back of the paper, and the fold’s numbers have to disappear. Finally, children cut along the last segment (continuous line). Two stars appear: a star hole in the sheet and a full star opening the cut paper. Figures \ref{fig:Starprimary} shows the students with both stars.

Figure \ref{fig:starpattern}(b) shows together the star and the lines that we fold so that one can see that these lines are symmetry axes of the star.

\item Activity 6: Geogebra construction. 

At middle school level, students investigate the steps of the construction of the star CP and reproduce it with Geogebra, see Figure \ref{fig:starpattern}(b) . This is very interesting because they have to construct the star (one idea is to start from a pentagon), to choose symmetry axes and to imagine what happens in each folding step. They renforce the knowledge about perpendicular bisector and bisectors and how they are constructed.

\begin{remark}

The scope of the use of this theorem could be extended at university level, since its proof is useful to understand algorithms complexity.
For other origami-based lessons for undergraduate students see \cite{Hul13}.
\end{remark}
\end{itemize}

\section{Results and discussion}
In this paper we presented a vertical project carry out principally in an Italian school, involving about 170 students from 6 to 12 years old. We propose a set of origami and ICT activities on a contemporary math result. The aim was to propose to students a series of problems that would allow them to deal with a mathematical subject from different points of view and with different approaches, based on active learning.
After the experimentation, following an exchange of considerations with students and teachers, we can say that the strengths of the project were as follows:
\begin{itemize}
\item[-] the topic is welcomed in all school level
\item[-] students work in class and at home and share and discuss their result (activities 1, 3, 4, 5); students work in group (activity 6).
\item[-] students can approach a problem with a tangible tool, they formalize and use ICT with awareness
\item[-] also linguistic skills are involved: students have to take care to different languages: origami  (CP and folding diagram reading), mathematics (specific terminology), English language.
\end{itemize}

Moreover, part of these activities were proposed for generic public, as a part of dissemination projects in Spain (at a session in the Museum of Mathematics of Catalunya) and in Italy (exhibition of a project at a school). In general is not so easy to propose maths activities to adults because of their fear of mathematics. But the topic and the way to propose it was well suited to an involvement of the adult public.
Here the comment of an adult:
\emph{“I would like to thank the teacher that during the exhibition showed to me how to cut a star with only one-cut. I proposed the same game to my friends in the evening and I looked great!”}
This comment shows also that we generate a second level of dissemination.

\section*{Acknowledgement(s)}

We want to thank students and teachers of the school “C. Gallazzi” and “C. Costamagna” in Busto Arsizio (Italy) for their support in this project.

First author was partially supported by the National Group of Algebraic and Geometric Structures and
their application (GNSAGA), Italy.


\end{document}